\documentclass[10pt]{article}

\usepackage{amsfonts,amsmath,theorem,latexsym,comment,euscript,amssymb}

\newcommand{\ah}{\"{a}}
\newcommand{\HH}{\EuScript{H}}
\newcommand{\Lb}{\ell}
\newcommand{\LL}{L^k}
\newcommand{\Lk}{\ell^{\otimes k}}
\newcommand{\Kk}{K^{(k)}}
\newcommand{\Qk}{Q^{(k)}}
\newcommand{\Qks}{Q^{(k)}_{\text{KS}}}
\newcommand{\bas}{\theta^{(k)}}
\newcommand{\tbas}{\widetilde{\theta}}
\newcommand{\lvol}{\epsilon_\omega}
\newcommand{\cd}{\varepsilon^{(k)}}
\newcommand{\cm}{\mu^{(k)}}
\newcommand{\cs}{\Phi^{(k)}}
\newcommand{\pks}{S^{(k)}}
\newcommand{\ncs}{\widetilde{\Phi}^{(k)}}
\newcommand{\xk}{x^{(k)}}
\newcommand{\Hk}{\EuScript{H}_k}
\newcommand{\pr}{\operatorname{proj}}

\newcommand{\Op}{\operatorname{Op}}
\newcommand{\Tr}{\operatorname{Tr}}

\newcommand{\lspan}{\operatorname{span}}
\newcommand{\rcs}{e^{(k)}}
\newcommand{\ta}{\psi^{(k)}}
\newcommand{\Ss}{\EuScript{S}}
\newcommand{\lefthook}{\,\text{\raisebox{0.44ex}{\rule{4pt}{0.3pt}}}\!\overset{\lrcorner}{}\,}
\newcommand{\dbar}{\overline{\partial}_k}
\newcommand{\sdbar}{\overline{\partial}\!\!\!/_k}
\newcommand{\ev}{\operatorname{ev}}
\newcommand{\real}{\operatorname{Re}}
\newcommand{\imag}{\operatorname{Im}}
\newcommand{\spinc}{Spin$^c$\ }
\newcommand{\zbar}{\bar{z}}
\newcommand{\wbar}{\bar{w}}
\newcommand{\medwedge}{\mbox{\fontsize{10pt}{0pt}\selectfont $\wedge$}}

\newcommand{\proj}[1]{\mathbb{P}#1}
\newcommand{\hprod}[2]{\langle #1 \vert #2 \rangle\,}
\newcommand{\hform}[2]{h(#1, #2)}
\newcommand{\ket}[1]{| #1 \rangle}
\newcommand{\bra}[1]{\langle #1 |}

\theorembodyfont{\rmfamily}
\newtheorem{thm}{Theorem}[section]
\newtheorem{definition}[thm]{Definition}
\newtheorem{prop}[thm]{Theorem}
\newtheorem{cor}[thm]{Corollary}
\newcommand{\proof}[1]{\smallskip\par\noindent\textbf{Proof: } #1 $\hfill\blacksquare$\par\smallskip}

\begin{document}

\begin{flushright}math.SG/0502026\end{flushright}
%\begin{flushright}Preprint\\ \today\end{flushright}
\vspace{7.5ex}
\begin{center}
 {\LARGE\rm Coherent States in Geometric Quantization} \\
 \vspace{7ex}
 {\large William D.\ Kirwin\footnote{E-mail address: {\tt  kirwin@math.usu.edu}}}\\
 \vspace{1.75ex}
 {\em Department of Mathematics and Statistics\\
 Utah State University, Logan, UT 84322-3900, USA}
\end{center}
\vspace{2ex}

\begin{abstract}  In this paper we study overcomplete systems of coherent states associated to compact integral symplectic manifolds by geometric quantization.  Our main goals are to give a systematic treatment of the construction of such systems and to collect some recent results.  We begin by recalling the basic constructions of geometric quantization in both the K\ah hler and non-K\ah hler cases.  We then study the reproducing kernels associated to the quantum Hilbert spaces and use them to define symplectic coherent states.  The rest of the paper is dedicated to the properties of symplectic coherent states and the corresponding Berezin-Toeplitz quantization.  Specifically, we study overcompleteness, symplectic analogues of the basic properties of Bargmann's weighted analytic function spaces, and the `maximally classical' behavior of symplectic coherent states.  We also find explicit formulas for symplectic coherent states on compact Riemann surfaces.\\

\end{abstract}

\begin{section}{Introduction}

Coherent states are ubiquitous in the mathematical physics literature. Yet there seems to be a lack of general theory in the context of geometric quantization.  This paper is an attempt to partially fill this gap.

We will define coherent states associated to an arbitrary integral compact symplectic manifold $(M,\omega)$ using the machinery of geometric quantization.  The metaplectic correction will not play a role in this construction and will be omitted for simplicity.  In the non-K\ah hler case there are at least two common methods of quantizing $M$: almost K\ah hler \cite{BU96,MM04} and \spinc \cite{Duis,MM02,M,TZ} quantization.  The definition of coherent states which we will describe is similar in both cases, although for technical reasons is somewhat simpler in the almost K\ah hler case.  In contrast to the quantization of a K\ah hler manifold, in both \spinc and almost K\ah hler quantization a quantum state does not necessarily have a nice holomorphic local form.  As a consequence, it is difficult to control the global behavior of the quantum states and, as we will see, the condition that $M$ is compact becomes essential.  This is not to say that the properties of the coherent states are different in the non-compact case -- it is simply not clear to the author how to proceed (see \cite{MM04} for some recent progress in this direction).

The definition of coherent states that we will make is semi-constructive.  It will depend on a choice of basis for the quantum Hilbert space arising from geometric quantization.  Since such a basis is not always available, finding an explicit form for the corresponding coherent states may be difficult.  On the other hand, the abstract approach taken here demonstrates that many of the traditional properties of coherent states follow from general considerations.  We will refer to the coherent states constructed here as symplectic coherent states in order to distinguish them from specific instances.

We should point out that symplectic coherent states are not, in general, of Perelomov-type \cite{Per}; i.e. they are not orbits of a fiducial vector under the action of a Lie group.  In some specific cases symplectic coherent states and Perelomov coherent states coincide -- for example in the quantizations of $\mathbb{C}$ and $S^2$ (see \S\ref{sec:examples}), where the methods of Perelomov yield a reproducing system in the quantum Hilbert space that arises from geometric quantization.  In fact, when this happens, the two construction must agree (Theorem \ref{prop:kuniq}).

The coherent state map we will study has appeared in a different (but in some ways equivalent) form in \cite{BU00} where it is used to prove a symplectic analogue of Kodaira's embedding theorem.  In \cite{BU00}, the coherent states lie in a different quantum Hilbert space and are associated to a circle bundle over $M.$  In \cite{MM04}, Ma-Marinescu analyze the semiclassical properties of generalized Bergman kernels.  This analysis leads to the notion of a peak section, which is related to symplectic coherent states (\S\ref{subsec:cs}). In \cite{Rawns77}, Rawnsley globalized the constructions of Perelomov \cite{Per} to make a general definition of a coherent state on K\ah hler manifolds.  When applied to the specific case of a principal $\mathbb{C}^\times$-bundle over $M$ the coherent states constructed here and the Rawnsley type coherent states are related (see \S\ref{subsec:rcs}).  The properties of Rawnsley-type coherent states, as well as the geometric interpretation of Berezin quantization which they provide, are studied in \cite{CGR1,CGR2}.  Symplectic coherent states are associated to $M$ itself, and can be associated to a non-K\ah hler symplectic manifold.

Symplectic coherent states generalize many of the systems constructed in the literature.  The most basic examples of coherent states are those associated to the complex plane -- the simplest K\ah hler manifold.  These states are known as Segal-Bargmann-Heisenberg-Weyl (or some permutation thereof) coherent states, or often more simply as canonical coherent states; see \cite{D80}, for example, where they are developed using projective representations of the symplectic group on the quantum Hilbert spaces of \S\ref{subsec:kahlerquant}.  Canonical coherent states are briefly described in \S\ref{sec:examples}.  \cite{KS} contains a survey on many traditional mathematical aspects of canonical coherent states,  the introduction of which also includes a discussion of what, in general, should be called a coherent state.

Following \cite{KS} (and to some extent popular opinion) we will define a system of coherent states to be a set $\{\ket{x}\in\HH\,\big\vert\,x\in M\}$ of quantum states in some quantum Hilbert space $\HH,$ parameterized by some set $M,$ such that:
\begin{enumerate}
\item the map $x\mapsto\ket{x}$ is smooth, and
\item the system is overcomplete; i.e.
\[\int_M \ket{x}\bra{x}\,d\mu(x)=\mathbf{1}_\HH.\]
\end{enumerate}
Physicists usually call property (2) completeness.  As we will see, the map $x\mapsto\ket{x}$ is actually antiholomorphic in a sense appropriate to non-K\ah hler manifolds.  The parameterizing set $M$ is generally, and for us will be, a classical phase space; i.e. an integral symplectic manifold.

We motivate our construction of coherent states by recalling some basic quantum mechanics.  The following observations are well known (\cite{KS} and \cite[Chap. 3]{Lands}).
Let us, for a moment, eschew definitions and rigor in order to see how to proceed.  The position space wave function representing a state $\ket{\psi}$ is $\psi(x)=\hprod{x}{\psi},$ where $\ket{x}$ is a coherent state localized at $x.$  The position space wave function of the coherent state $\ket{x}$ is then
\[K(x,y):=K_x(y):=\hprod{y}{x}.\]
We can use this to rewrite the equation $\hprod{x}{\psi}=\psi(x)$ as
\begin{equation}
\label{eqn:rep1}
\int \overline{K(x,y)}\psi(y)\,dy=\psi(x).
\end{equation}
A function $K$ that satisfies \eqref{eqn:rep1} for some space of functions is called a reproducing kernel for that space.  Reading the above discussion backwards, we see that a coherent state can be defined in terms of a reproducing kernel.  We will use this approach and define symplectic coherent states in terms of reproducing kernels for the quantum Hilbert spaces arising from geometric quantization, sometimes called generalized Bergman kernels; properties (1) and (2) will then follow.  This construction is well known for the K\ah hler quantization of $\mathbb{C}$ and yields the Bergman reproducing kernel, which in turn yields the canonical coherent states.  The asymptotics of generalized Bergman kernels are studied in \cite{BU96,BU00,DLM,MM02,MM04} and will be important when we consider the semiclassical limit.

\smallskip

\noindent\textbf{Example} A reproducing kernel for the space $\EuScript{S}(\mathbb{R})$ of Schwartz functions on $\mathbb{R}$ is the Dirac distribution $\delta(x-y).$ \hfill$\blacksquare$

\smallskip

Symplectic coherent states associated to the Poincar\'{e} disc, and hence via Riemann uniformization to compact Riemann surfaces of genus $g\geq2,$ were used in \cite{KL1,KL2} to study the semiclassical limit of the deformation quantizations of these surfaces.  We will give explicit formulas for these coherent states in \S\ref{sec:examples}.

The rest of the paper is organized as follows:  in \S\ref{sec:background} we review the geometric quantization of a K\ah hler manifold $(M,\omega)$ and two generalizations to the non-K\ah hler case known as the almost K\ah hler and \spinc quantizations of $(M,\omega).$  In \S\ref{sec:cs} we define the reproducing kernel for the quantum Hilbert space associated to $M$ by geometric quantization and use it to define symplectic coherent states.  We also describe symplectic analogues of some analytic function space results of \cite{Barg61}, and the relationship between Rawnsley-type and symplectic coherent states.  In \S\ref{sec:btq} we discuss the overcompleteness relation and the coherent state quantization induced by the symplectic coherent states.  In \S\ref{sec:cb} we show that symplectic coherent states are the most classical quantum states and consider the semiclassical limit.  Finally, in \S\ref{sec:examples} we apply the constructions of \S\ref{sec:cs} to compact Riemann surfaces.

\end{section}

\begin{section}{Background and Notation}
\label{sec:background}

  \begin{subsection}{Prequantization}

Throughout we assume that $(M,\omega)$ is an integral compact symplectic manifold; i.e. $\left[\frac{\omega}{2\pi}\right]$ is in the image of the map $H^2(M;\mathbb{Z})\rightarrow H^2_{\text{DR}}(M).$  The basic object of geometric quantization is an Hermitian line bundle $\pi:\Lb\rightarrow M$ with compatible connection $\nabla$ with curvature $-i\omega,$ known as the prequantum line bundle.  The existence of $\Lb$ is guaranteed by the integrality of $\omega;$ in fact the Chern character of $\Lk$ is $ch(\Lk)=\exp(k\left[\frac{\omega}{2\pi}\right]).$  For detailed accounts of geometric quantization see \cite{GS,Woo91}.

We denote by $h:\overline{\Lb}_x\otimes \Lb_x\rightarrow\mathbb{C}$ the Hermitian structure on $\Lb.$  We follow the physics convention that the first term is conjugate linear.  All tensor products will be taken over $\mathbb{C}.$  The norm of $q\in L_x$ is $|q|^2=\hform{q}{q}.$  $h$ induces an Hermitian structure on $\Gamma(\Lb):$ for $s_1,s_2\in\Gamma(\Lb)$
\[\hprod{s_1}{s_2}=\int_M\hform{s_1(x)}{s_2(x)}\lvol(x),\]
where
\[\lvol=\left(\frac{1}{2\pi}\right)^n\frac{\omega^{\wedge n}}{n!}\]
 is the Liouville volume form on $M.$  We have included recurring factors of $2\pi$ in the Liouville form.  This will simplify some formulas later on, but will also have the effect that our formulas differ slightly from some of those in the literature.  The norm of $s\in\Gamma(\Lb)$ is $\|s\|^2=\hprod{s}{s}.$

We will be occasionally interested in the semiclassical $k=1/\hbar\rightarrow\infty$ limit of $\Lk.$  The structures $(h,\hprod{\cdot}{\cdot},\nabla)$ on $\Lb$ induce corresponding structures on $\Lk$ which we will denote by the same symbols.  The curvature of $\Lk$ is $-ik\omega.$

The program of geometric quantization associates to $(M,\omega)$ a Hilbert space $\HH$ and a map $Q:C^\infty(M)\rightarrow\Op(\HH).$  To begin, we define the prequantum Hilbert space $\Hk^0$ to be the $L^2$ completion, with respect to the Liouville measure, of the set of square integrable sections of $\Lk:$
\[\Hk^0=\overline{\{s\in\Gamma(\Lk)\,\big\vert\,\|s\|<\infty\}}.\]

The Kostant-Souriau quantization of the Poisson-Lie algebra $C^\infty(M)$ is the map
\[f\in C^\infty(M)\mapsto \Qks(f)=-\frac{i}{k}\nabla_{X_f}+f \in\Op(\Hk^0)\]
where $X_f$ is the Hamiltonian vector field defined by
\[X_f\lefthook\omega=df.\]
In \S\ref{sec:btq} we will recall an alternate quantization of $C^\infty(M).$

As is well known, $\Hk^0$ is too large for the purposes of quantization \cite[Ch 9]{Woo91}.  If $(M,\omega)$ is K\ah hler there is a standard method of choosing a subspace $\Hk\subset\Hk^0.$  In the non-K\ah hler case, there are (at least) two reasonable methods: almost K\ah hler quantization \cite{BU96,MM04} and \spinc quantization \cite{Duis,MM02,M,TZ}.  We will review these three constructions in the next three sections.

  \end{subsection}

  \begin{subsection}{K\ah hler Quantization}
\label{subsec:kahlerquant}

If $(M,\omega,J)$ is a K\ah hler manifold with complex structure $J,$ there is a natural method of reducing the prequantum Hilbert space $\Hk^0.$  The complexified tangent bundle of $M$ decomposes into $\pm i$ eigenspaces of $J:$
\begin{equation}
\label{eqn:tmdecomp}
TM_\mathbb{C}=TM_J^{(1,0)}\oplus TM_J^{(0,1)}.
\end{equation}
A section $s\in\Gamma(\Lk)$ is said to be polarized if it is tangent to the K\ah hler polarization $TM_J^{(1,0)};$ i.e. if $\nabla_X s=0$ for each $X\in\Gamma(TM^{(0,1)}).$

The quantum Hilbert space is defined to be the $L^2$ closure of the set of polarized sections of $\Lk:$
\[\Hk=\overline{\{s\in\Gamma(\Lk)\,\big\vert\,\|s\|<\infty, s\text{ polarized}\}}.\]

This quantum Hilbert space can be described in terms of a Dirac-type operator \cite[Chap. 3]{BGV}.  The decomposition \eqref{eqn:tmdecomp} induces a decomposition
\begin{equation}
\label{eqn:extdecomp}
\Lambda^*(T^*M) =\underset{p,q=0}{\overset{n}{\oplus}}\Lambda^{p,q}(T^*M) =\underset{p,q=0}{\overset{n}{\oplus}} \Lambda^p(T^*M_J^{(1,0)})\otimes\Lambda^q(T^*M_J^{(0,1)}).
\end{equation}
Let $\dbar:\Omega^{p,q}(M,\Lk)\rightarrow\Omega^{p,q+1}(M,\Lk)$ denote the Dolbeault operator twisted by $\Lk.$  Hodge's theorem says that $\ker(\dbar+\dbar^*)^2$ is isomorphic to the sheaf cohomology space $H^*(M,\mathcal{O}(\Lk)).$  Kodaira's vanishing theorem then tells us that for $k$ sufficiently large, $H^q(M,\mathcal{O}(\Lk))=0$ for $q>0,$ and hence that $\Hk=\ker\dbar\big\vert_{\Gamma(\Lk)}.$  The dimension of $\Hk$ can be computed with the Riemann-Roch-Hirzebruch theorem:
\begin{equation}
\label{eqn:RR}
d_k:=\dim\Hk=RR(M,\Lk)=\int_M ch(\Lk)Td(TM_J^{(1,0)}).
\end{equation}
In particular, since we assume $M$ is compact, $d_k<\infty.$

  \end{subsection}

  \begin{subsection}{Almost K\ah hler Quantization}

We suppose now that $(M,\omega)$ is an integral compact symplectic manifold, not necessarily K\ah hler.  Every such manifold admits an $\omega$-compatible almost complex structure $J,$ and any two such choices are homotopic.  The complexified tangent bundle again decomposes as in \eqref{eqn:tmdecomp}.  If $J$ is not integrable (i.e. $M$ is not a complex manifold), then there may be no sections $s\in\Gamma(\Lk)$ which are tangent to $TM_J^{(1,0)}$.  In this case, more work is required to define a quantum Hilbert space.

In this section, we describe one such method, known as almost K\ah hler quantization, which was introduced in \cite{BU96} based on results in \cite{GU} and further studied in \cite{DLM,MM02,MM04}.  The idea is to replace $(\dbar+\dbar^*)^2,$ which does not exist if $J$ is not integrable, with the rescaled Laplacian $\Delta_k:=\Delta-nk,$ where $\Delta$ is the Laplacian for the metric $g=\omega(\cdot,J\cdot).$  In the K\ah hler case, these two quantities are equated by the Bochner-Kodaira formula: $\Delta_k=2(\dbar+\dbar^*)^2.$  The main result of \cite{GU} is (in a slightly sharpened form due to \cite{BU96,MM04}):

\begin{thm} Given an integral symplectic manifold $(M,\omega)$ with $\omega$-compatible almost complex structure, there exists a constant C and a positive constant $a$ such that for $k$ sufficiently large,
\begin{enumerate}
\item the first $d_k$ eigenvalues of $\Delta_k$ (in nondecreasing order) lie in the interval $(-a,a),$ and
\item the remaining eigenvalues lie to the right of $nk+C,$
\end{enumerate}
where $d_k$ is the Riemann-Roch number as in \eqref{eqn:RR}.
\end{thm}

Motivated by this, the quantum Hilbert space is defined to be
\[\Hk:=\lspan_\mathbb{C}\{\bas_1,\dots,\bas_{d_k}\}\]
where $\bas_j$ is an eigenfunction of $\Delta_k$ with eigenvalue $\lambda_j$ and $\lambda_1\leq\lambda_2\leq\cdots\leq\lambda_j\leq\cdots.$  Note that just as in the K\ah hler case, the dimension is given by the Riemann-Roch number $d_k.$   This quantization has excellent semiclassical ($\hbar=1/k\rightarrow0$) properties \cite{BMS,BU96,MM02,MM04}.  It also has the advantage that the quantum Hilbert space consists of sections of $\Lk.$  \spinc quantization, described in the next section, does not have either of these properties, although the \spinc quantum states are true zero-modes of a Dirac-type operator.

In analogy with the K\ah hler case, we will call the elements of $\Hk$ polarized sections.  Since $M$ is compact and the polarized sections are smooth and the quantum Hilbert space is finite dimensional, $\Hk$ is a closed subspace of $L^2(M,\Lk).$

  \end{subsection}

  \begin{subsection}{\spinc Quantization}

The idea of \spinc quantization is to find a suitable generalization of $\dbar$ for the non-K\ah hler case.  We will briefly review the relevant details here.  See \cite[App D]{LM} for a more complete account of the \spinc bundle, and \cite{Duis,MM02,M,TZ}  for \spinc quantization.

Let $J$ be an $\omega$-compatible almost complex structure on $(M,\omega)$ so that $g=\omega(\cdot,J\cdot)$ is a Riemannian metric on $M.$  The \spinc bundle associated to the data $(M,\omega,J)$ is defined as $S(M):=\Lambda^{0,*}(T^*M)$ according to the decomposition \eqref{eqn:extdecomp}.  There is a Dirac-type operator $\sdbar:\Omega^{0,*}(M,\Lk)\rightarrow\Omega^{0,*}(M,\Lk)$ that decomposes into $({\sdbar})_+:\Omega^{0,\text{even}}(M,\Lk) \rightarrow\Omega^{0,\text{odd}}(M,\Lk)$ and $(\sdbar)_-:\Omega^{0,\text{odd}}(M,\Lk) \rightarrow\Omega^{0,\text{even}}(M,\Lk).$

The quantum Hilbert space associated to the data $(M,\omega,J)$ is the virtual vector space
\[\Hk:=\ker(\sdbar)_+\ominus\ker(\sdbar)_-.\]
The dimension is again given by the Riemann-Roch number $d_k$ \eqref{eqn:RR}.  There is a \spinc analogue of the Kodaira vanishing theorem \cite{BU96,MM02} which insures that $\Hk$ is an honest vector space for $k$ sufficiently large.

The metric $g$ on $M$ and the Hermitian structure $h$ on $\Lk$ combine to give an Hermitian structure, also denoted by $h,$ on $S(M).$  Although the zero-modes of $\sdbar$ are not sections of $\Lk$ since they have higher degree components, their norms are asymptotically concentrated on the zero degree part; i.e. there exists a constant $C>0$ such that for $k$ sufficiently large $\|s_+\|\leq C k^{-1/2}\|s_0\|,$ for each $s\in\Hk$ and where $s=s_0+s_+$ denotes the decomposition of $s$ into zero and higher degree components \cite{BU96,MM02}.

Just as in the K\ah hler and almost K\ah hler cases, we will refer to the elements of $\Hk$ as polarized sections.  Also, since $M$ is compact and polarized sections are again smooth, and the space of them is finite dimensional, $\Hk$ is a closed subspace of $L^2(M,S(M)).$

We will assume throughout that $k$ is chosen sufficiently large to ensure the validity of the relevant vanishing/existence theorem.

  \end{subsection}

\end{section}

\begin{section}{Coherent States}
\label{sec:cs}

In this section we will construct coherent states associated to an integral symplectic manifold $(M,\omega).$

  \begin{subsection}{Reproducing Kernels}
  \label{subsec:repkern}

For the K\ah hler, almost K\ah hler and \spinc quantizations of a compact symplectic manifold $(M,\omega)$ the quantum Hilbert space is finite dimensional with dimension given by the Riemann-Roch formula \eqref{eqn:RR}.  We will see below that since $d_k<\infty$ there exists a reproducing kernel for the quantum Hilbert space.  In the K\ah hler case, the existence of a reproducing kernel may also be established by trivializing $\Lk$ -- polarized sections are locally holomorphic and standard methods from complex analysis can be used.  In the non-K\ah hler cases no such nice local form is known to the author, and the compactness assumption and resulting finite dimensionality become essential.

For two vector bundles $\pi_j:E_j\rightarrow M_j,\ j=1,2$ we define $E_1\boxtimes E_2:=\pi_1^*E_1\otimes\pi_2^*E_2\rightarrow M_1\times M_2.$  If a vector bundle $E$ has an Hermitian structure $h$, we will identify $\overline{E}\simeq E^*$.  For $u,v\in E_x$ we define $\overline{u}\cdot v:=\hform{u}{v}.$  Similarly, we identify $\overline{u}\otimes u=\hform{u}{u}.$   Moreover, if $L$ is a line bundle, we will identify $\overline{v}\otimes w=\hform{v}{w}$ for $\overline{v}\in\overline{L}_x,\ w\in L_x$ so that $\overline{L}\otimes L\simeq\mathbb{C}.$  Combining these definitions, we see that $u\otimes\overline{v}\cdot w=\hform{v}{w}u.$  These conventions agree since $u\otimes\overline{v}\otimes w=\hform{v}{w}u=\hform{v}{u}w=u\otimes\overline{v}\cdot w.$  Finally, we define $\overline{\overline{u}\otimes v}:=u\otimes\overline{v}.$

Many of the results in this section hold for all three methods of quantization.  To unify notation and avoid repetition we define
\[
\LL:=\begin{cases}
\Lk & \text{ for Kahler and almost Kahler quantization}\\
\Lambda^{0,*}(T^*M)\otimes\Lk & \text{ for \spinc quantization}.
\end{cases}
\]

Let $\{\bas_j\}_{j=1}^{d_k}$ be a unitary basis for $\HH_k$.

\begin{definition}
The reproducing kernel $\Kk\in\Gamma(\overline{\LL}\boxtimes\LL)$ is the section
\[\Kk(x,y):=\sum_{j=1}^{d_k}\overline{\bas_j(x)}\otimes\bas_j(y).\]
\end{definition}

$\Kk$ is also known as a generalized Bergman kernel.  Note that $\Kk$ does not depend on the choice of unitary basis.

In the K\ah hler and \spinc quantization schemes, the quantum Hilbert space is the kernel of a Dirac-type operator, and the reproducing kernel is the large $t$ limit of the associated heat kernel (see \S\ref{subsec:quant}).  Although it is not a function on $M\times M$, the reproducing kernel $\Kk$ has many of the same properties enjoyed by reproducing kernels for analytic function spaces.

\begin{prop}
\label{prop:kuniq}
$\Kk$ is the unique polarized section of $\overline{\LL}\boxtimes\LL$ such that
\[\int_M \overline{\Kk(x,y)}\cdot s(y)\,\lvol(y) = s(x)\quad \forall s\in\HH_k.\]
\end{prop}

\proof{Suppose there are two reproducing kernels.  Their difference evaluated against an arbitrary section $s\in\Hk$ must be zero.  Hence this difference must be in $\overline{\Hk}\otimes(\Hk)^\bot$.  But since both kernels are polarized, the difference is in $\overline{\Hk}\otimes\Hk$ and is therefore zero.}

The restriction of $\Kk$ to the diagonal is a smooth function.

\begin{definition}
\label{def:cd}
The coherent density is the smooth function $\cd\in C^\infty(M)$ defined by
\[\cd(x):=\Kk(x,x)=\sum_{j=1}^{d_k}|\bas_j(x)|^2.\]
Since $M$ is compact and $\cd$ is smooth and nonnegative, we may define a measure on $M$ by $\cm=\cd\lvol$ which we will call the coherent measure.
\end{definition}

Since $\HH_k\subseteq L^2(M,\LL)$ is a closed subspace, there is a projection $\Pi_k:L^2(M,\LL)\rightarrow\HH_k.$  To find the Schwartz kernel of this projection, we need the following observation (which follows from the facts that $L^2(\mathbb{R}^{2n},\mathbb{C})$ is separable and that since $M$ is compact it has a finite cover by open sets which are diffeomorphic to subsets of $\mathbb{R}^{2n}$).

\begin{prop} $L^2(M,\LL)$ is a separable Hilbert space.
\end{prop}

\begin{prop}
The Schwartz kernel of $\Pi_k$ is the reproducing kernel $\Kk.$
\end{prop}

\proof{We need to show that
\[\left(\Pi_k s\right)(x)=\int_M\overline{\Kk(x,y)}\cdot s(y)\lvol(y)\qquad\forall s\in L^2(M,\LL).\]
Since $\Pi_k$ is a projection, it is uniquely characterized by $\text{im }\Pi_k=\Hk$ and $\Pi_k^2=\Pi_k=\mathbf{1}_{\Hk}.$

Let $\{\bas_j\}_{j=1}^\infty$ be a unitary basis for $L^2(M,\LL)$ such that $\lspan_\mathbb{C}\{\bas_j\}_{j=1}^{d_k}=\Hk.$  Then for each $s\in L^2(M,\LL)$ there exists $\{s^j\in\mathbb{C}\}_{j=1}^\infty$ such that $\|s-\sum_{j=1}^Ns^j\bas_j\|^2\rightarrow0$ as $N\rightarrow\infty.$  We then have
\begin{equation}
\label{eqn:interchange1}
\int_M\overline{\Kk(x,y)}\cdot s(y)\lvol(y)=\sum_{j=1}^{d_k}\bas_j(x)\int_M\sum_{l=1}^\infty s^l \hform{\bas_j(y)}{\bas_l(y)}\lvol(y).
\end{equation}
Using H\"{o}lder's inequality we see that
\[\sum_{l=1}^\infty |s^l|^2 \int_M |\hform{\bas_j(y)}{\bas_l(y)}|\lvol(y)\leq\sum_{l=1}^\infty |s^l|^2\|\bas_j(y)\|^2\|\bas_l(y)\|^2=\|s\|^2.\]
Hence, the integrand is absolutely integrable and we may interchange the integral and sum in \eqref{eqn:interchange1} to obtain
\[\int_M\overline{\Kk(x,y)}\cdot s(y)\lvol(y)=\sum_{j=1}^{d_k}\bas_j(x)\in\Hk\]
as desired.

Moreover, for $s=\sum_{j=1}^{d_k}s^j\bas_j(x)\in\Hk$, we easily obtain $\Pi_k^2 s=\Pi_k s$ in terms of $\Kk$ since all of the relevant sums are finite.  We conclude that $\Kk$ is the Schwartz kernel of $\Pi_k$ as desired.}
\end{subsection}

\begin{subsection}{Coherent States}
    \label{subsec:cs}

We now define the coherent states associated to an integral compact symplectic manifold $(M,\omega).$

\begin{definition}
\label{def:cs}
The coherent state localized at $x\in M$ is
\[\cs_x:=\Kk(x,\cdot) = \sum_{j=1}^{d_k}\overline{\bas_j(x)}\otimes\bas_j.\]
\end{definition}

In order to distinguish these coherent states from others, we will sometimes refer to them as symplectic coherent states.  Observe that $\cs$ depends smoothly and antiholomorphically on $x$ in the generalized sense: a section is holomorphic on a symplectic manifold if it is polarized.

Since $\cs_x\in\overline{\LL_x}\otimes\Hk,$ it is necessary to investigate how $\cs_x$ should be interpreted as a quantum state.  Consider the case of almost K\ah hler quantization.  If we trivialize $\Lk_x$ with a unit, $\cs_x$ becomes a well-defined state in $\Hk$ via the identification $1\otimes\bas_j\simeq\bas_j.$  The different unit trivializations of $\Lk_x$ are parameterized by $U(1)$ which means that $\cs_x$ is a well-defined quantum state up to a phase -- the usual situation in quantum mechanics.  On the other hand, quantum states are most properly regarded as rays in the projective Hilbert space $\proj{\Hk}.$  It follows from the above discussion that the map $x\in M\mapsto\mathbb{C}\cdot\cs_x\in\proj{\Hk}$ is well-defined and smooth.  If it is not possible to find a global unit section of $\Lk$ then there is no smooth lift of $\cs$ to $\Hk.$

In \cite{MM04}, this map is shown to be asymptotically symplectic (as $k\rightarrow\infty$), asymptotically isometric with respect to the metric $g=\omega\circ J,$ and, for $k$ sufficiently large, an embedding (see \cite{BU00} for similar results).

We now return our attention to the general case of almost K\ah hler or \spinc quantization.  It is sometimes convenient to work with normalized states.  In terms of coherent states, Definition \ref{def:cd} reads

\begin{prop}
\label{thm:normcs}
$\|\cs_x\|^2=\cd(x).$
\end{prop}

This is the reason $\cd$ is called the coherent density.

Let $M_k:=\{x\in M\,\big\vert\,\cd(x)\neq0\}.$  By Corollary \ref{cor:zeroes}, $M_k$ is the complement of the base locus of $\Hk.$  For $x\in M_k$ we will denote the normalized coherent state localized at $x$ by
\[\ncs_x:=\ket{\xk}:=\frac{\cs_x}{\sqrt{\cd(x)}}.\]
For $x\not\in M_k$ we define $\ncs_x:=\ket{\xk}:=0.$
We can now state the reproducing property concisely:
\begin{gather*}
\hprod{\cs_x}{s}=s(x)\\
\sqrt{\cd(x)}\hprod{\xk}{s}=s(x)\qquad \forall x\in M,\forall s\in\Hk.
\end{gather*}
For $x\not\in M_k$ the above is justified by Corollary \ref{cor:zeroes}.

If $x\in M_k$ then $\mathbb{C}\cdot\cs_x=\mathbb{C}\cdot\ket{\xk}\in\proj{\Hk}$ so that we may define the projection $\pr_{\ket{\xk}}:\Hk\rightarrow\Hk$ by
\begin{equation}
\label{eqn:prj}
s\mapsto(\cd(x))^{-1}\ket{\cs_x}\hprod{\cs_x}{s} =\ket{\xk}\hprod{\xk}{s}.
\end{equation}
We will also write $\pr_{\mathbb{C}\cdot\cs_x} =\ket{\xk}\bra{\xk} =(\cd(x))^{-1}\ket{\cs_x}\bra{\cs_x}$ for this projection.

The following is a generalization of a result in \cite{Barg61} and is the basic reason why the quantum Hilbert space behaves in many ways like a weighted analytic function space, even on non-K\ah hler manifolds.  We will further develop this analogy in \S\ref{subsec:overcomp}.

\begin{prop}
\label{prop:basicineq}
For each polarized section $s\in\HH_k,$
\[|s(x)|\leq\|s\|\sqrt{\cd(x)}.\]
\end{prop}

\proof{Let $s\in\Hk.$  We use the reproducing property of the coherent states to write, for $x\in M_k,$
\[|s(x)|^2=\hform{s(x)}{s(x)}=\hform{\hprod{\cs_x}{s}}{\hprod{\cs_x}{s}} =|\hprod{\cs_x}{s}|^2.\]
The result then follows from the Cauchy-Schwartz inequality and Theorem \ref{thm:normcs}.

If $x\not\in M_k$ then $\cd(x)=0$ implies $\bas_j(x)=$ for all $j.$  Since $s$ is a linear combination of $\bas_j,$ this implies $s(x)=0.$ }

This Theorem has two useful corollaries.  The proofs are simple and left to the reader.

\begin{cor}
\label{cor:evalcont}
The evaluation map $\ev_x:s\in\HH_k\mapsto s(x)\in\LL_x$ is continuous.
\end{cor}

We can also prove Corollary \ref{cor:evalcont} directly -- it follows from the facts that $d_k<\infty$ and that the sets $\{|\bas_j(x)|\}_{j=1}^{d_k}$ and $\{\|\bas_j\|\}_{j=1}^{d_k}$ are bounded, which in turn follow from our assumption that $M$ is compact.

In the K\ah hler case, the compactness assumption on $M$ can be lifted.  The existence of the reproducing kernel, as well as Theorem \ref{prop:basicineq}, can then be deduced from Jensen's formula \cite[p324]{Lang}: if $f$ is holomorphic on the closed disc of radius $R$ and $f(0)\neq0,$ and the zeroes of $f$ in the open disc, ordered by increasing moduli and repeated according to multiplicity, are $z_1,...,z_N,$ then
\[|f(0)|\leq \frac{\|f\|_R}{R^N}|z_1\cdots z_N|.\]

Unfortunately, the author is unaware of nice local forms of polarized sections in the \spinc and almost K\ah hler quantizations of a non-K\ah hler manifold.  This makes the assumption that $M$ is compact, or more precisely that $d_k<\infty$, essential.  Of course, it is also not clear whether the spectrum of the rescaled Laplacian has the requisite structure to define the almost K\ah hler quantization of $M$ in the non-compact case.

On the other hand, most of the results of this paper hold for any choice of finite dimensional subspace of the prequantum Hilbert space $\Hk^0.$  The primary advantages of almost K\ah hler and \spinc quantization are that they provide canonical methods for choosing such a subspace and that they provide enough structure to ensure a meaningful semiclassical limit (\emph{c.f.} \eqref{eqn:pksasymp}).

The next corollary justifies our definition of the normalized coherent state $\ket{\xk}$ in the case that \\ $\cd(x)=0.$

\begin{cor}
\label{cor:zeroes}
$\cd(x)=0$ if and only if $s(x)=0$ for each polarized section $s.$
\end{cor}

In \S\ref{subsec:ta}, \S\ref{sec:cb}, and \S\ref{sec:examples} we will be interested in the semiclassical limit of the symplectic coherent states.  It will be useful to express $\cs_x$ in terms of the peak sections of Ma-Marinescu \cite{MM04}, defined as follows:  The Kodaira map $\Psi^{(k)}:M_k\rightarrow \proj{\Hk^*}$, which sends $x\in M_k$ to the hyperplane $\{s\in\Hk\,\big\vert\,s(x)=0\}$ of sections which vanish at $x$, is base point free for large enough $k.$  Construct an unitary basis $\{\bas_1,\dots,\bas_{d_k-1},\pks_x\}$ such that $\bas_j(x)=0$ for $1\leq j\leq d_k-1.$  Then $\pks_x,$ called a peak section, is a unit norm generator of the orthogonal complement of $\Psi^{(k)}(x).$  Observe that
\begin{equation}
  \label{eqn:cspks}
    \cs_x(y)=\Kk(x,y)=\overline{\pks_x(x)}\otimes\pks_x(y)
\end{equation}
and also that $\cd(x)=|\pks_x(x)|^2.$  Moreover,
\begin{equation}
  \label{eqn:pksnorm}
    \int_M |\pks_x(y)|^2\lvol(y)=1.
\end{equation}

  \end{subsection}

  \begin{subsection}{Rawnsley-Type Coherent States}
\label{subsec:rcs}

The coherent states defined in \cite{Rawns77} for compact K\ah hler manifolds are a generalization of Bargmann's principal vectors \cite{Barg61} to spaces of holomorphic sections.  We will describe their relation to symplectic coherent states in this section.  Due to Corollary \ref{cor:evalcont}, we are able to construct Rawnsley-type coherent states on any compact integral symplectic manifold.  In this section, we will consider only the almost K\ah hler quantization of $M$ so that the prequantum bundle is $\Lk.$

By Corollary \ref{cor:evalcont}, for each $q\in \Lb_x$ we get a continuous map $\delta_q:\Hk\rightarrow\mathbb{C}$ by composing the evaluation $\ev_x$ with the trivialization $s(x)=\delta_q(s)q^{\otimes k}.$  By the Riesz representation theorem, there exists $\rcs_q\in\Hk$ such that
\[\hprod{\rcs_q}{s}q^{\otimes k}=\delta_q(s)q^{\otimes k} = s(\pi(q))\qquad \forall s\in\Hk.\]
Observe that
\[\rcs_{cq}=\overline{c}^{-k}\rcs_q\]
for $0\neq c\in\mathbb{C}.$  We will refer to the section $\rcs_q$ as a Rawnsley-type coherent state.  We can express the reproducing kernel, and therefore the symplectic coherent states, in terms of $\rcs_q:$

\begin{prop}
\label{prop:kernrcs}
Let $x\in M$ and $q\in\Lb_x.$  Then
\[\Kk(x,y)=\overline{q}^{\otimes k}\otimes\rcs_q(y).\]
Equivalently, $\cs_x=\overline{q}^{\otimes k}\otimes\rcs_q.$
\end{prop}

\proof{In terms of the unitary basis $\{\bas_j\}_{j=1}^{d_k}$ for $\Hk$ we have
\begin{equation*}
s(x)=\hprod{\cs_x}{s} = \sum_{j=1}^{d_k}\hprod{\bas_j}{s}\bas_j(x)
= \sum_{j=1}^{d_k}\hprod{\bas_j}{s}\tbas_j^q(x)q^{\otimes k}
=\hprod{\rcs_q}{s}q^{\otimes k}
\end{equation*}
where $\tbas_j^q$ is the trivialization of $\bas_j$ determined by a local section with value $q^{\otimes k}$ at $x.$  Therefore
\begin{equation}
\label{eqn:Rcsexpansion}
\rcs_q = \sum_j\overline{\tbas_j^q(\pi(q))}\bas_j.
\end{equation}
Hence we obtain,
\begin{equation*}
\overline{q}^{\otimes k}\otimes \rcs_q = \cs_{\pi(q)}.
\end{equation*}}

Let $s_0:M\rightarrow \Lb^\times.$  In \cite{Rawns77}, Rawnsley defines a function
\[\eta(x):=\hprod{e^{(1)}_{s_0(x)}}{e^{(1)}_{s_0(x)}}|s_0(x)|^2.\]
It is easy to check that this function is independent of $s_0.$  This function was also studied for K\ah hler $M$ in \cite{CGR1,CGR2} and in the almost K\ah hler case in \cite{BU00}.  A short calculation using the previous Theorem yields:

\begin{cor} $\eta=\varepsilon^{(1)}.$
\end{cor}

  \end{subsection}

  \begin{subsection}{Transition Amplitudes}
    \label{subsec:ta}

In this section we define the 2-point transition amplitude for symplectic coherent states and show that it can be interpreted as a probability density on $M.$

\begin{definition} The 2-point function, or transition amplitude, is
\[\ta(x,y):=|\hprod{\xk}{y^{(k)}}|^2\in C^\infty(M\times M).\]
\end{definition}

In terms of the reproducing kernel, the 2-point function is
\[\ta(x,y)=\frac{\Kk(y,x)\cdot\Kk(x,y)}{\cd(x)\cd(y)}.\]
As expected, $\ta(x,x)=1.$  The Cauchy-Schwartz inequality
\[|\hprod{\cs_x}{\cs_y}|^2\leq\|\cs_x\|^2\|\cs_y\|^2\]
implies $\ta(x,y)\in[0,1].$  Since the map $x\mapsto\mathbb{C}\cdot\cs_x$ is an embedding for $k$ sufficiently large \cite{MM04}, $x\neq y$ implies $\cs_x\neq\cs_y$ so that $\ta(x,y)=1$ if and only if $x=y.$

\begin{prop}
\label{prop:taprobdensity} For each $x\in M_k,$
$\ta(x,y)$ is a probability density on $M$ with respect to the coherent measure $\cm(y).$
\end{prop}

\proof{For each $x$ with $\cd(x)\neq0,$
\begin{equation}
\label{eqn:2ptprob}
\int_M \ta(x,y)\,d\cm(y) = \frac{1}{\cd(x)}\int_M \Kk(y,x)\cdot\Kk(x,y)\lvol(y) = 1.
\end{equation}}

In \cite{CGR2}, Cahen-Gutt-Rawnsley define a 2-point function for K\ah hler $M$ in terms of Rawnsley-type coherent states:
\[\psi'(x,y) = \frac{|\hprod{e_q^{(1)}}{e_{q'}^{(1)}}|^2} {\|e_q^{(1)}\|^2\|e_{q'}^{(1)}\|^2}\]
where $x=\pi(q)$ and $y=\pi(q').$  By Theorem \ref{prop:kernrcs} we see that $\psi'=\psi^{(1)}.$  Moreover, if the quantization is regular (i.e. $\cd(x)=const$ for all $k$) then it follows from \cite[Prop. 2 and eq 1.7]{CGR2} that $\ta=(\psi^{(1)})^k.$

The transition amplitude can be expressed in terms of peak sections by a simple calculation using equation \eqref{eqn:cspks}:

\begin{prop}
\label{prop:tapks}
$\ta(x,y)\cm(y)=|\pks_x(y)|^2$
\end{prop}

Theorem \ref{prop:taprobdensity} is therefore equivalent to \eqref{eqn:pksnorm}.

Finally, we note here that the association to each $\hbar=1/k,\ k\in\mathbb{Z}_+$ of $\Hk,$ $\mathbb{C}\cdot\cs_x$ and $\cm$ defines a pure state quantization of the integral symplectic manifold $(M,\omega)$ (see \cite[p113]{Lands} for the definitions).

  \end{subsection}

\end{section}

\begin{section}{Berezin-Toeplitz Quantization}
\label{sec:btq}

In this section we study the Berezin quantization \cite{Ber74} induced by the coherent state map $\cs.$  This method of quantization is studied in detail in the context of analytic function spaces in \cite{Lands}.  The extension to Hilbert spaces of sections of the prequantum bundle can be described in terms of symplectic coherent states.

  \begin{subsection}{Overcompleteness and Characteristic Sets}
    \label{subsec:overcomp}

In this section we consider the most important property of coherent states: overcompleteness.

\begin{definition}
  A system of coherent states $\{\ket{x}\in\Hk\,\big\vert\,x\in M\}$ is overcomplete with respect to a measure $\mu$ if
\begin{enumerate}
  \item $\hprod{x}{y}\neq 0$ for all $x,y\in M$ with $\ket{x},\ket{y}\neq 0,$ and
  \item $\int_M \ket{x}\bra{x}\,d\mu(x)=\mathbf{1}_{\Hk}.$
\end{enumerate}
\end{definition}

\begin{thm}
\label{thm:overcomp}
The system of symplectic coherent states $\{\ket{\xk}\,\big\vert\,x\in M\}$ defined in \S\ref{sec:cs} is overcomplete with respect to the coherent measure $\cm.$  In particular,
\begin{equation}
\label{eqn:overcomp}
\int_M\ket{\xk}\bra{\xk}\,d\cm(x)=\mathbf{1}_{\Hk}.
\end{equation}
\end{thm}

\proof{We compute, for every $s_1,s_2\in\Hk,$
\begin{gather*}
\bra{s_1} \int_M \ket{\xk}\bra{\xk}\,d\cm(x) \ket{s_2}
= \int_M \hprod{s_1}{\cs_x}\hprod{\cs_x}{s_2}\,\lvol(x) \\
= \int_M \hform{s_1(x)}{s_2(x)}\,\lvol(x) = \hprod{s_1}{s_2}.
\end{gather*}
so that $\int_M \ket{\xk}\bra{\xk}\,d\cm = \mathbf{1}_{\Hk}$ as desired.}

\begin{cor}
There exist points $x_1,\dots,x_{d_k}\in M$ such that $\{\ket{\xk_1},\dots,\ket{\xk_{d_k}}\}$ is a basis for $\Hk.$
\end{cor}

\proof{Let $x_1\in M_k$ and set $\Ss_1=\{\ket{\xk_1}\}.$  If $d_k=1$ then $\Ss$ is a basis for $\Hk.$  Suppose $d_k>r\geq1$ and let $\Ss_r=\{\ket{\xk_1},\dots,\ket{\xk_r}\}$ be a set of linearly independent vectors in $\Hk.$  Since $d_k>r,$ there is some vector $\ket{\psi}\not\in\lspan_\mathbb{C} \Ss_r.$  Suppose for every $x\in M$ that $\ket{\xk}\in\lspan_\mathbb{C} \Ss_r.$  Then
\[\int_M \ket{\xk}\hprod{\xk}{\psi}d\cm(x)\in\lspan_\mathbb{C} S_r.\]
which implies
\[\int_M \ket{\xk}\hprod{\xk}{\psi}d\cm(x)\not=\ket{\psi}.\]
This contradicts Theorem \ref{thm:overcomp}.  Hence, there is some $x\in M$ such that $\ket{\xk}\not\in\lspan_\mathbb{C} \Ss_r.$  Let $x_{r+1}=x.$  Then $\Ss_{r+1}$ is a linearly independent set in $\Hk.$  We continue inductively.  Since $d_k<\infty,$  the process must stop, and the resulting set is the required linearly independent set.}

This corollary motivates the following definition, which is a generalization of the characteristic point sets introduced in \cite{Barg61}.

\begin{definition}
A set $\Ss\subseteq M$ is characteristic if for every $s\in\Hk,$
\[s\big\vert_\Ss=0 \text{ implies } s=0.\]
\end{definition}

\begin{prop} If $\Ss\subseteq M$ is characteristic, then $\{\ket{\xk}\,\big\vert\,x\in\Ss\}$ is complete.
\end{prop}

\proof{If $s(x)=0$ for all $x\in\Ss$ implies $s=0,$ then $\hprod{\xk}{s}=0$ for all $x\in\Ss$ implies $s=0.$  Hence, the only vector orthogonal to $\{\ket{\xk}\,\big\vert\,x\in\Ss\}$ is 0, which means $\Ss$ is complete.}

  \end{subsection}

  \begin{subsection}{Quantization}
    \label{subsec:quant}

In this section we study the quantizing map $Q:C^\infty(M)\rightarrow\Op(\Hk)$ resulting from the overcompleteness relation \eqref{eqn:overcomp}.  Applying Berezin's method of quantization \cite{Ber74}, we have

\begin{definition}
The Berezin quantization $\Qk(f)$ of $f\in C^\infty(M)$ is the operator
\[\Qk(f):=\int_M f(x)\ket{\xk}\bra{\xk}\,d\cm(x).\]
\end{definition}

In fact, $\Qk(f)$ converges for $f\in L^\infty(M).$  Some basic theorems about Berezin's method of quantization of analytic function spaces apply in this case;  see for example \cite[Thm 1.3.5]{Lands}.
There is another way to describe the Berezin quantization of $f.$  For each $s\in L^2(M,\LL)$ we have
\[\left(\Pi_k s\right)(x)=\int_M \overline{\Kk(x,y)}\cdot s(y)\lvol(y) = \hprod{\cs_x}{s}.\]
Therefore
\begin{align*}
(\Qk(f)s)(x) &= \left(\int_M f(y)\ket{y^{(k)}}\hprod{y^{(k)}}{s}\,d\cm(y)\right)(x) \\
&= \int_M \overline{\Kk(x,y)}\cdot (f(y)(\Pi_k s)(y))\lvol(y) =\left(\Pi_k\circ M_f\circ\Pi_k s\right)(x)
\end{align*}
where $M_f$ denotes the multiplication operator.  In this form, $\Qk(f)$ is known as the Toeplitz quantization of $f.$  The Berezin-Toeplitz quantization and Kostant-Souriau quantizations of $f$ are related by Tuynman's formula \cite{Tuy}:
\[\Pi_k\circ\Qk_\text{KS}(f)\circ\Pi_k =\Qk(f-\frac{1}{2k}\Delta f)\]
where $\Delta$ is the Laplacian associated to the metric $g=\omega(\cdot,J\cdot).$  See \cite{BdMG} for the theory of generalized Toeplitz operators, and \cite{BMS} for an analysis of the semiclassical properties of $\Qk$.

We can recast Berezin's covariant symbol \cite{Ber74} in terms of the symplectic coherent states.

\begin{definition}
\label{def:covsymb}
The covariant symbol $\widehat{A}\in C^\infty(M)$ associated to the operator $A\in\Op(\Hk)$ is
\[\widehat{A}(x):=\bra{\xk}A\ket{\xk},\]
where $A\ket{\xk}:=\sum_{j=1}^{d_k}\overline{\bas_j(x)}\otimes A\bas_j.$
\end{definition}

A consequence of Theorem \ref{prop:kernrcs} is that Definition \ref{def:covsymb} agrees with the covariant symbol defined in \cite{CGR1} for K\ah hler $M$ using Rawnsley-type coherent states.  A standard result involving Berezin's covariant symbol is true in our case as well (see \cite{CGR2} for an analogous computation with Rawnsley-type coherent states):

\begin{prop}
\label{prop:trace}
$\Tr A=\int_M\widehat{A}(x)\,d\cm(x).$\end{prop}

We conclude this section by pointing out the relationship, in the \spinc and K\ah hler cases, between symplectic coherent states and the heat kernel of the appropriate Dirac-type operator.  See \cite[Chap. 3]{BGV} for a detailed analysis of the heat kernel, some properties of which we will use below.

The heat kernel $\Kk_t\in C^\infty(\mathbb{R}_+\times M\times M,\overline{\LL}\boxtimes\LL)$ of the Laplacian associated to $\dbar$ (or $\sdbar$ in the \spinc case) admits an expansion
\begin{equation}
\label{eqn:heatkernexp}
\Kk_t(x,y)=\sum_{j=0}^\infty e^{-\lambda_j t}\ \overline{\bas_j}(x)\otimes\bas_j(y)
\end{equation}
where $0\leq\lambda_1\leq\lambda_2\leq\dots\rightarrow\infty$ are the eigenvalues of the Laplacian with corresponding eigenmodes $\bas_j\in\Gamma(\LL).$  Moreover,
\[\left|\sum_{j=d_k+1}^\infty e^{-t \lambda_j }\overline{\bas_j}\otimes\bas_j\right|\leq C e^{-\lambda_{d_k+1}t}\]
for some constant $C>0$ \cite[Prop 2.37]{BGV}.  Hence, the large time limit of the heat kernel is a symplectic coherent state:
\[\lim_{t\rightarrow\infty}\Kk_t(x,y)=\cs_x(y).\]

In the almost K\ah hler case, although the heat kernel has an expansion of the form \eqref{eqn:heatkernexp}, the low lying eigenvalues of the polarized states are not necessarily zero, and so the large time limit of the heat kernel is not directly related to the symplectic coherent states.

Finally, observe that $\cm=\underset{t\rightarrow\infty}{\lim}\Tr\Kk_t$ and so Theorem \ref{prop:trace}, applied to the identity operator, yields the familiar index formula $d_k=\int_M \underset{t\rightarrow\infty}{\lim} \Tr\Kk_t.$

  \end{subsection}

\end{section}

\begin{section}{Classical and Semiclassical Behavior}
\label{sec:cb}

\begin{subsection}{Classical Behavior for Finite $k$}
\label{subsec:classical}
In this section we show that the coherent states defined in \S\ref{sec:cs} are the quantum states that behave most classically: they are maximally peaked and evolve classically.

Consider for a moment the almost K\ah hler quantization of $M$ so that the prequantum bundle $\LL$ is the line bundle $\Lk.$  In this case, the coherent state $\cs_x$ is the projection of the Dirac distribution onto $\Hk.$  To see why, let $x\in M$ and trivialize $\Lk$ over an open set $U$ containing $x$ with a unit section $s_0.$  Let $\tilde{\delta}^{(k)}_x$ denote the Dirac distribution on $U$ centered at $x$ and define
\[\delta_x^{(k)}(y) =\begin{cases}\tilde{\delta}_x^{(k)}(y)\overline{s_0(x)}\otimes s_0(y) & \text{ for }y\in U \\
0 & \text{ otherwise}
\end{cases}
\]
Then
\begin{equation*}
\hprod{\delta_x^{(k)}}{s}=s(x)\quad \forall s\in C^1(M,\Lk).
\end{equation*}
$\delta_x^{(k)}$ does not depend on our choice of $s_0$.  We now see that $\cs_x$ is the projection onto $\Hk$ of $\delta_x^{(k)}$ since
\[\left(\Pi_k\delta_x^{(k)}\right)(y) =\int_M\overline{\Kk(y,z)}\cdot\delta_x^{(k)}(z)\lvol(z)=\cs_x(y).\]

We next observe that symplectic coherent states are maximally peaked quantum states.  The following result holds for coherent states arising from K\ah hler, almost K\ah hler and \spinc quantization.

\begin{thm} $\cs_x$ maximizes $|s(x)|^2$ over all $s\in\overline{\Lk_x}\otimes\Hk$ with $\|s\|^2=\cd(x).$
\end{thm}

\proof{As in Theorem \ref{prop:basicineq} we write
\[|s(x)|^2=|\hprod{\cs_x}{s}|^2.\]
This is minimized when we have equality in the Cauchy-Schwartz inequality, which occurs when $s$ is proportional to $\cs_x.$}

In the almost K\ah hler case, we can say more: $\cs_x$ evolves classically.  Suppose $f\in C^\infty(M)$ is such that the Hamiltonian vector field $X_f$ is complete.  The flow of $X_f$ induces a Hamiltonian evolution of sections in the quantum Hilbert space as follows \cite[\S 8.4]{Woo91}:

We can lift the Hamiltonian vector field $X_f$ to a vector field $V_f$ on $T(\Lk);$ that is, there exists a unique vector field on $\Lk$ defined by $\pi_* V_f=X_f$ and $\frac{1}{k}V_f\lefthook\alpha=\frac{1}{k}V_f\lefthook\overline{\alpha} =f\circ\pi,$ where $\alpha$ is the connection 1-form on the complement of the zero section of $\Lk.$  If the fiber coordinate is $z=re^{i\phi}$ in a local trivialization, then
\[V_f=X_f+k L\frac{\partial}{\partial \phi}\]
where $L=X_f\lefthook\tau+f$ is the Lagrangian associated to $f$ by a local symplectic potential $\tau.$  Locally, $T\Lk\simeq TM\times\mathbb{C}$ and we have identified $X_f\in TM$ and $\frac{\partial}{\partial \phi}\in T\mathbb{C}$ with the corresponding vector fields in $TM\times\mathbb{C}.$  The flow $\xi_t$ of $V_f$ is fiber preserving and projects to the flow $\rho_t$ of $X_f.$  Moreover, $\xi_t$ induces a linear pull-back action $\widehat{\rho}_t:\Gamma(\Lk)\rightarrow\Gamma(\Lk)$ by
\begin{equation}
\label{eqn:vflow}
\xi_t(\widehat{\rho}_t s(x))=s(\rho_t x).
\end{equation}
In fact, this action is infinitesimally generated by the Kostant-Souriau quantization of $f:$
\[\frac{d}{dt}\widehat{\rho}_t=ik\Qk_\text{KS}(f)\widehat{\rho}_t.\]
This is one of the motivations for the Kostant-Souriau quantization $\Qk_\text{KS}$ of $f.$  The restriction of $\widehat{\rho}_t$ to $\Hk^0$ is a 1-parameter unitary group.

Extending the action of  $\xi_t$ to $\overline{\Lk}\boxtimes\Lk$ in the obvious way, we have:

\begin{prop}
For $f\in C^\infty(M)$ such that $X_f$ is complete,
\begin{equation}
\label{eqn:classicalevolution}
\xi_t\cs_x(y)=\cs_{\rho_t x}(\rho_t y);
\end{equation}
i.e. the symplectic coherent states evolve classically.  Equivalently,
\[\widehat{\rho}_t\cs_x=\cs_x.\]
\end{prop}

\proof{Since $\widehat{\rho}_t$ is a unitary endomorphism of $\Hk^0,$ we have \begin{align*}
\widehat{\rho}_t\cs_x(y)&=\sum_{j=1}^{d_k} \overline{\widehat{\rho}_t\bas_j(x)}\otimes\widehat{\rho}_t\bas_j(y) \\
&=\sum_{j=1}^{d_k}\overline{\left(\bas_j(x)\right)} \,^t\overline{\widehat{\rho}_t} \otimes \widehat{\rho}_t\bas_j(y)\\
&=\cs_x(y).
\end{align*}}

\end{subsection}

\begin{subsection}{The Semiclassical Limit}

The asymptotic analysis of generalized Bergman kernels by Ma-Marinescu reveals the semiclassical behavior of peak sections \cite[eq 3.24]{MM04}.  Let $\{r_k\}$ be a sequence of real numbers with $r_k\rightarrow0$ and $\sqrt{k}\ r_k\rightarrow\infty$ as $k\rightarrow\infty.$  Denote by $B(x,r)$ the open geodesic ball of radius $r$ centered at $x\in M.$  Then
\begin{equation}
  \label{eqn:pksasymp}
    \int_{B(x,r_k)}|\pks_x(y)|^2\lvol(y) =1-O(1/k),\qquad k\rightarrow\infty.
\end{equation}
Comparing this with \eqref{eqn:pksnorm} we see that the peak section $\pks_x$ is asymptotically concentrated about $x.$  In terms of the transition amplitude, \eqref{eqn:pksasymp} is
\begin{equation}
  \label{eqn:taasymp}
    \int_{B(x,r_k)}\ta(x,y)\,\cm(y)=1-O(1/k),\qquad k\rightarrow\infty.
\end{equation}
Combining this with Theorem \ref{prop:taprobdensity}, we have:

\begin{prop}
\label{prop:sclimit}
If $f\in C^1(M)$ then
\[\lim_{k\rightarrow\infty}\int_M f(y)\ta(x,y)\cm(y) = f(x);\]
that is, $\lim_{k\rightarrow\infty}\ta(x,y)\cm(y)=\delta_x(y).$
\end{prop}

\proof{Let $\{r_k\}$ be a sequence of positive real numbers with $r_k\rightarrow0$ and $\sqrt{k}\ r_k\rightarrow\infty.$  By Theorem \ref{prop:tapks} we have, for each $x\in M_k,$
\begin{equation}
  \label{eqn:sclimitint}
  \begin{split}
    \big\vert \int_M(f(y)-f(x)) & \ta(x,y)\,d\cm(y)\big\vert \\
    &\leq \int_{B(x,r_k)} |f(y)-f(x)|\,|\pks_x(y)|^2\lvol(y) + \int_{M\setminus B(x,r_k)} |f(y)-f(x)|\,|\pks_x(y)|^2\lvol(y).
  \end{split}
\end{equation}
The first integral on the right hand side of \eqref{eqn:sclimitint} goes to zero as $k\rightarrow\infty$ because of \eqref{eqn:pksasymp} and the fact that $f$ is continuous.  The second integral on the right hand side of \eqref{eqn:sclimitint} goes to zero since $f(x)-f(y)$ is bounded (specifically as a function of $y$) and equations \eqref{eqn:pksasymp} and \eqref{eqn:pksnorm} imply that the peak sections go to zero outside the ball $B(x,r_k).$
}

Of course, there is a physical reason to expect this behavior.  The coherent state localized at $x$ should be the quantum state most concentrated about $x.$  In the semiclassical limit, we expect to recover the classical picture -- in particular the classical state most concentrated about $x$ is the Dirac distribution at $x.$

\end{subsection}

\end{section}

\begin{section}{Examples}
\label{sec:examples}

\begin{subsection}{The Complex Plane}
This example is well known \cite{KS,Woo91}.  We will take $z=\frac{1}{\sqrt{2}}(x+i y),$ use the standard symplectic form $\omega=idz\medwedge d\zbar,$ and trivialize the prequantum line bundle (globally since $\mathbb{C}$ is contractible) with the symplectic potential $\tau=\frac{i}{2}(z d\zbar-\zbar d z).$  The space $\Hk=\{f(z)e^{-k|z|^2/2}\}$ of polarized sections of $\Lk$ relative to the standard complex structure can be identified with a weighted Bargmann space \cite{Barg61}.  A unitary basis for $\Hk$ is $\{\sqrt{\frac{k}{j!}}\,z^j e^{-k|z|^2/2}\}_{j\in\mathbb{N}}.$  The reproducing kernel is the usual Bergman kernel
\[\Kk(w,z)=k\sum_{j=0}^\infty \frac{(\wbar z)^j}{j!}e^{-k|z|^2/2-k|w|^2/2} = k e^{\wbar z-k|z|^2/2-k|w|^2/2}.\]
The coherent density is $\cd=k$ and the 2-point function is $e^{-k|z-w|^2}.$  In this case, it is easy to see that the semiclassical limit yields the expected results:
\[\lim_{k\rightarrow\infty} \ta(w,z)\cm(z) =\lim_{k\rightarrow\infty}\left(\frac{k}{2\pi}\right)e^{-k|w-z|^2}\omega =\delta(w-z)\omega.\]
\end{subsection}

\begin{subsection}{The 2-Sphere}
Coherent states on $S^2$ are constructed by Perelomov in \cite{Per} using Lie group techniques.  As we will see, the construction of \S\ref{sec:cs} yields the same results without using any group structure.  The correspondence between the two methods is due to Theorem \ref{prop:kuniq} and the fact that Perelomov's coherent states are reproducing.

We trivialize $S^2\simeq\mathbb{C}P^1=\{[z_0,z_1]\}/\mathbb{C}$ over the open set $U_0=\{[z_0,z_1]\,\big\vert\,z_0\neq0\}.$  Define a local coordinate $z=z_1/z_0.$  The Fubini-Study symplectic form on $U_0$ is
\[\omega=\frac{idz\medwedge d\zbar}{(1+|z|^2)^2}.\]
Trivializing $\Lk$ with the symplectic potential $\tau=(1+|z|^2)^{-1} i\zbar dz,$
the Hermitian form on $\Lk$ is given by $h(p,q)=(1+|z|^2)^{-k}\overline{p}q.$
We then have the following unitary basis for $\Hk:$
\[\left\{\sqrt{(k+1)\binom{k}{j}}\ z^j\,\big\vert\,0\leq j\leq k\right\}.\]
The coherent state localized at $w$ is therefore
\[\cs_w(z)=(k+1)\sum_{j=0}^k\binom{k}{j}(\wbar z)^j =(k+1)(1+\wbar z)^k.\]
The corresponding coherent density is $\cd=k+1;$  one must take care to include the extra factors arising from the Hermitian structure when evaluating $\cs_w(w).$  In this case, Theorem \ref{prop:sclimit} becomes
\begin{equation}
\lim_{k\rightarrow\infty}\ta(w,z)\cm(z)=\lim_{k\rightarrow\infty}\frac{k+1}{2\pi} \left[\frac{(1+\wbar z)(1+w\zbar)}{(1+|z|^2)(1+|w|^2)}\right]^k \frac{i dz\medwedge d\zbar}{(1+|z|^2)^2}=\delta(w-z).
\end{equation}
\end{subsection}

\begin{subsection}{The 2-Torus}
For $\lambda=\lambda_1+i\lambda_2\in\mathbb{C}$ with $\operatorname{Im}\lambda>0,$ let $T^2(\lambda)=\mathbb{C}/\{m+n\lambda\,\big\vert\,m,n\in\mathbb{Z}\}.$  $\lambda$ is known as the modulus of the torus.  The standard symplectic form $\omega=2\pi i\lambda_2^{-1}dz\medwedge d\zbar$  on $\mathbb{C},$  normalized to be integral on $T^2(\lambda),$  descends to a symplectic form on $T^2(\lambda).$  The prequantum line bundle $\Lk$ can be lifted to a line bundle over $\mathbb{C}$ (since $\mathbb{C}$ is the universal cover of $T^2(\lambda)$).  The resulting line bundle can be globally trivialized.  Hence we will identify sections of $\Lk$ with appropriately pseudoperiodic functions on $\mathbb{C}.$

If we trivialize with the symplectic potential $\tau=i\pi\lambda_2^{-1}(z d\zbar-\zbar dz)$ then the Hermitian form is given by $h(p,q)=\overline{p}q.$  A unitary basis for the quantum Hilbert space can be given in terms of $\vartheta$-functions.  Let
\begin{align*}
\vartheta_j(\lambda;z)&=\sum_{n\in\mathbb{Z}} e^{i\lambda\pi(kn^2+2jn)+2\pi i\sqrt{2}(j+k n)z}, \\
\psi^{(k)}_j(z,\zbar)&=e^{k\pi z(z-\zbar)/\lambda_2} \vartheta_j(\lambda;z), \text{ and} \\
N_{k,j}=\|\psi^{(k)}_j\|^2&=\frac{1}{\sqrt{2k\lambda_2}}\,e^{2\pi j^2\lambda_2/k}.
\end{align*}
A unitary basis for $\Hk$ is $\{N_{k,j}^{-1/2}\psi^{(k)}_j(z,\zbar)\}_{j=0}^{k-1}.$ From this we construct the coherent state localized at $w=\frac{1}{\sqrt{2}}(w_1+i w_2):$
\[\cs_w(z)=\sum_{j=0}^{k-1}\sqrt{2k\lambda_2}\, e^{-2\pi j^2\lambda_2/k} e^{\sqrt{2}\pi i k(zy-\wbar w_2)/\lambda_2} \overline{\vartheta_j(\lambda;w)}\vartheta_j(\lambda;z).\]
The coherent density is
\[\cd(z)=\sqrt{2k\lambda_2} e^{-2\pi k y^2/\lambda_2} \sum_{j=0}^{k-1}e^{-2\pi j^2\lambda_2/k} |\vartheta_j(\lambda;z)|^2.\]
The semiclassical limit (Theorem \ref{prop:sclimit}) yields the identity
\begin{align*}
\delta(w-z)=\lim_{k\rightarrow\infty}\ta(w,z)\cm(z)&=\lim_{k\rightarrow\infty} \sqrt{2k\lambda_2}\,e^{-2\pi k w_2^2/\lambda_2} \left(\sum_{j=0}^{k-1}e^{-2\pi j^2\lambda_2/k}|\vartheta_j(\lambda;z)|^2\right)^{-1} \\ &\qquad\qquad \cdot \sum_{j,l=0}^{k-1} e^{-2\pi(j^2+l^2)\lambda_2/k} \overline{\vartheta_j(\lambda;w)\vartheta_l(\lambda,z)} \vartheta_j(\lambda;z)\vartheta_l(\lambda;w)\,\frac{2\pi i}{\lambda_2}dz\medwedge d\zbar.
\end{align*}

\end{subsection}

\begin{subsection}{Higher Genus Riemann Surfaces}
We will construct coherent states on a compact Riemann surface $\Sigma_g$ of genus $g\geq2$ by uniformizing $\Sigma_g$ as the quotient of the complex upper half plane $\mathfrak{H}=\{z\in\mathbb{C}\,\big\vert\,\operatorname{Im}z>0\}$ by a Fuschian group $\Gamma$ (see \cite{AMV,Jost} for details).  The coherent states in this section correspond to those of Klimek-Lesniewski \cite[eq 4.5]{KL2}.

Let $\Gamma<PSL(2,\mathbb{Z})$ be a Fuschian group.  $PSL(2,\mathbb{Z}),$ and hence $\Gamma,$ acts on $\mathfrak{H}$ by fractional linear transformations.  The space $\Sigma_g:=\Gamma\backslash\mathfrak{H}$ is a compact manifold if and only if $\Gamma$ is a hyperbolic group, which we will henceforth assume.  The K\ah hler form $\omega=i(\operatorname{Im}z)^{-2} dz\medwedge d\zbar$ descends to a symplectic form on $\Sigma_g,$ as do the complex structure and K\ah hler metric.

For $\gamma=\begin{pmatrix}a & b\\c & d\end{pmatrix}\in\Gamma$ and $z\in\mathfrak{H}$ define $j(\gamma,z):=cz+d.$  An automorphic form of weight $k$ relative to $\Gamma$ is a function $f:\mathfrak{H}\rightarrow\mathbb{C}$ such that $f(\gamma z)=j(\gamma,z)^kf(z).$  The space $\mathcal{A}_k(\mathfrak{H})$ of automorphic forms on $\mathfrak{H}$ of weight $k$ relative to $\Gamma$ has an Hermitian product
\[\hprod{f}{g}_k=\int_\mathfrak{H}\overline{f(z)}g(z) (\operatorname{Im}z)^k\frac{\omega}{2\pi}.\]

The bundle $\Lb:=T^*\Sigma_g^{(1,0)}$ is a prequantum bundle for $\Sigma_g$ since its curvature is $-i\omega.$  The quantum Hilbert space $\Gamma(\Sigma_g,\mathcal{O}(T^*\Sigma_g^{(1,0)}))$  is isomorphic to $\mathcal{A}_2(\Sigma_g).$  Sections of $\Lk$ correspond to automorphic forms of weight $2k$ relative to $\Gamma$ restricted to $\Sigma_g.$  The Hermitian form descends to $\Sigma_g$ and is known as the Weil-Petersson inner product.

The reproducing kernel for $\mathcal{A}_{2k}(\mathfrak{H})$ is known (see for example \cite{Mi}) and descends to $\Sigma_g$ via a Poincar\'{e} series.  The resulting coherent state localized at $w$ is
\[\cs_w(z)=\frac{k-1}{2}\sum_{\gamma\in\Gamma}\left(\frac{2i}{\gamma z-\wbar}\right)^k j(\gamma,z)^{-k}.\]
The associated coherent density is
\[\cd(z)=(k-1)\sum_{\{\gamma,\gamma^{-1}\}\subset\Gamma} \real\left[\frac{2i}{(\gamma z-\zbar)j(\gamma,z)}\right]^{2k}.\]
Finally, the semiclassical limit of Propostion \ref{prop:sclimit}:
\begin{align*}
  \delta(w-z)=\lim_{k\rightarrow\infty}&\ta(w,z)\cm(z)\\ &=\lim_{k\rightarrow\infty}\frac{k-1}{2\pi}\,4^{k-1}(\imag z)^{k-2}(\imag w)^k\left(\sum_{\{\gamma,\gamma^{-1}\}\subset\Gamma} \real\left[\frac{2i}{(\gamma z-\zbar)j(\gamma,z)}\right]^{2k}\right)^{-1} \\ &\qquad\qquad\cdot\sum_{\gamma,\gamma'\in\Gamma} (|\gamma z-\wbar|\,|j(\gamma,z)|)^{-2k}\,i dz\medwedge d\zbar.
\end{align*}

\end{subsection}

\end{section}

\noindent\textbf{Acknowledgments.}  The author is grateful to Siye Wu for many valuable discussions and suggestions.  He would also like to thank Xiaonan Ma for bringing to his attention the papers \cite{DLM,MM02,MM04} as well as for several helpful comments.

\bibliographystyle{amsplain}
\bibliography{cs}

\end{document}